\documentclass[11pt]{article}
\usepackage{amsmath,amssymb}
\newtheorem{theo}{Theorem}[section]
\newtheorem{prop}[theo]{Proposition}
\newtheorem{lemma}[theo]{Lemma}

\newtheorem{conj}[theo]{Conjecture}

\newcommand{\ignore}[1]{}
\def\square{\vrule height6pt width7pt depth1pt}
\def\endpf{\hfill\square\bigskip}

\begin{document}
\title{The number of edge-disjoint transitive triples in a tournament}
\author{Raphael Yuster
\thanks{
e-mail: raphy@research.haifa.ac.il \qquad
World Wide Web: http:$\backslash\backslash$research.haifa.ac.il$\backslash$\symbol{126}raphy}
\\ Department of Mathematics\\ University of
Haifa at Oranim\\ Tivon 36006, Israel}

\date{} 

\maketitle
\setcounter{page}{1}
\begin{abstract}
We prove that a tournament with $n$ vertices has more than $0.13n^2(1+o(1))$ edge-disjoint transitive triples.
We also prove some results on the existence of large packings of $k$-vertex transitive
tournaments in an $n$-vertex tournament. Our proofs combine probabilistic arguments and some powerful packing results
due to Wilson and to Frankl and R\"odl.
\end{abstract}

\section{Introduction}
All graphs and digraphs considered here are finite and have no loops or multiple edges.
For the standard terminology used the reader is referred to \cite{Bo}.
A {\em tournament} on $n$ vertices is an orientation of $K_n$.
Thus, for every two distinct vertices $x$ and $y$, either $(x,y)$ or $(y,x)$ is an edge,
but not both.

Let $TT_k$ denote the unique transitive tournament on $k$ vertices. $TT_3$
is also called a {\em transitive triple} as it consists of some triple $\{(x,y), (x,z), (y,z)\}$.
A {\em $TT_k$-packing} of a directed graph $D$ is a set of edge-disjoint copies of $TT_k$
subgraphs of $D$. The {\em $TT_k$-packing number} of $D$, denoted $P_k(D)$,
is the maximum size of a $TT_k$-packing of $D$.
The $TT_3$-packing number of $D_n$, the complete digraph with $n$ vertices and $n(n-1)$ edges,
has been extensively studied. See, e.g., \cite{Ga,PhLi,Sk}.

In this paper we consider only $TT_k$-packings of tournaments.
Let $f_k(n)$ denote the minimum possible value of $P_k(T_n)$, where $T_n$ ranges over all possible
$n$-vertex tournaments. For simplicity, put $f(n)=f_3(n)$ and $P(T_n)=P_3(T_n)$.
Trivially, $P_k(T_n) \leq n(n-1)/(k(k-1))$, and in particular $f(n) \leq n(n-1)/6 < 0.167n^2(1+o(1))$.
In fact, it is not difficult to show that $f(n) \leq \lceil n(n-1)/6-n/3 \rceil$ (see Section 4
for this and also for a general way to construct an upper bound for $f_k(n)$).
We conjecture the following:
\begin{conj}
\label{c1}
$f(n)= \lceil n(n-1)/6-n/3 \rceil$.
\end{conj}
This conjecture was verified for all $n \leq 8$.
Our main result is the following lower bound for $f(n)$.
\begin{theo}
\label{t1}
$f(n) > 0.13n^2(1+o(1))$.
\end{theo}
We prove Theorem \ref{t1} in Section 2. In section 3 we show that if $T_n$ is
the random tournament on $n$ vertices then $P(T_n) \geq \frac{1}{6}n^2(1-o(1))$ almost surely.
In fact, we show that $P_k(T_n) \geq \frac{1}{k(k-1)}n^2(1-o(1))$ almost surely.
The final section contains some concluding remarks.

\section{Proof of the main result}
From here on we assume that the vertex set of a tournament with $k$
vertices is $[k]=\{1,\ldots,k\}$. Let $T_k$ be any $k$-vertex tournament.
For $v \in [k]$, let $d^+(v)$ denote the out-degree of $v$ in $T_k$.
Let $a(T_k)$ denote the total number of transitive triples in $T_k$,
and let $t(T_k)$ denote the total number of directed triangles in $T_k$.
Clearly, $a(T_k)+t(T_k)={k \choose 3}$. We shall also make use of the obvious inequality,
which follows from the fact that in a transitive triple there is one source and one sink.
\begin{equation}
\label{e1}
a(T_k) = \sum_{i=1}^k \frac{1}{2}\left({{d^+(v)} \choose 2} + {{k-1-d^+(v)} \choose 2}\right)
\geq \frac{k(k-1)(k-3)}{8}.
\end{equation}
In the proof of Theorem \ref{t1} we need the following (special) case
of Wilson's Theorem \cite{Wi}.
\begin{lemma}
\label{l21}
There exists a positive integer $N$ such that for all $n > N$, if $n \equiv 1 \bmod 42$
then $K_n$ decomposes into ${n \choose 2}/21$ edge-disjoint copies of $K_7$.
\endpf
\end{lemma}
The next lemma quantifies the fact that if $t(T_7)$ is relatively small then $P(T_7)$ is relatively large.
\begin{lemma}
\label{l22}
If $t(T_7) \leq 4$ then $P(T_7)=7$. If $t(T_7) \leq 11$ then $P(T_7) \geq 6$.
If $t(T_7) \geq 12$ then $P(T_7) \geq 5$.
\end{lemma}
{\bf Proof}\,
Clearly, the expected number of directed triangles in a random Steiner triple system
of $T_7$ is $7 \cdot \frac{t(T_7)}{t(T_7)+a(T_7)}= \frac{7}{35}t(T_7)$.
Hence, if $t(T_7) < 5$ then this expectation is less than $1$. Thus, there is a Steiner
triple system with no directed triangle. Namely, $P(T_7)=7$ in this case.
Similarly, by \ref{e1}, we always have $a(T_7) \geq 21$ and so $t(T_7) \leq 14$.
Therefore, the expectation above is always at most $14 \cdot \frac{7}{35} \leq 2.8$.
Thus, there is always a Steiner triple system with at most two directed triangles.
Namely, $P(T_7) \geq 5$ always.

We remain with the case where $t(T_7) \leq 11$. Notice that we may assume $t(T_7)=11$
or $t(T_7)=10$ since otherwise the above expectation argument yields $P(T_7) \geq 6$.
Assume first that $t(T_7)=11$. Hence $a(T_7)=24$ and by (\ref{e1}) the only possible scores
(sorted out-degree sequence) of such a $T_7$ are $(4,4,4,3,2,2,2)$,
$(5,3,3,3,3,2,2)$ and $(4,4,3,3,3,3,1)$. The last two scores are complementary (namely,
reversing the edges of a $T_7$ with one of these scores yields a tournament with the other score)
and the first score is self-complementary. Hence, one needs only to check the first two scores.

There are precisely 18 non-isomorphic tournaments with the score $(4,4,4,3,2,2,2)$,
and each can be checked to have at least 6 edge-disjoint transitive triples.
A convenient way to enumerate these 18 non-isomorphic tournaments is as follows.
Let $A_i$ be the set of vertices with out-degree $i$, $i=2,3,4$. $|A_2|=|A_4|=3$, $|A_3|=1$.
First case: The subgraph induced by $A_2$ is a directed triangle and the subgraph
induced by $A_4$ is also a directed triangle. There are four non-isomorphic tournaments
with this restriction. Second case: The subgraph induced by $A_2$ is a directed triangle
and the subgraph induced by $A_4$ is a transitive triple. There are four non-isomorphic tournaments
with this restriction. Third case: The subgraph induced by $A_2$ is a transitive triple
and the subgraph induced by $A_4$ is a directed triangle. There are four non-isomorphic tournaments
with this restriction. Fourth case: Both $A_2$ and $A_4$ induce a transitive triple. There are six
non-isomorphic tournaments with this restriction. Altogether there are $4+4+4+6=18$ possibilities.

There are precisely 15 non-isomorphic tournaments with the score $(5,3,3,3,3,2,2)$,
and each can be checked to have at least 6 edge-disjoint transitive triples.
A convenient way to enumerate these 18 non-isomorphic tournaments is as follows.
Let $A_i$ be the set of three vertices with out-degree $i$, $i=2,3,5$. $A_2=\{a,b\}$
$A_5=\{c\}$, $A_3=\{d,e,f,g\}$. We may assume the edge inside $A_2$ is $(a,b)$.
First case: $(a,c)$ is an edge. There is a unique tournament with this restriction.
Second case: $(c,a)$, $(a,d)$ and $(d,c)$ are edges. There are four non-isomorphic
tournaments. Third case: $(c,a)$, $(a,d)$ $(c,d)$ and $(b,d)$ are edges.
There are three non-isomorphic tournaments. Fourth case: $(c,a)$, $(a,d)$ $(c,d)$ and $(d,b)$
are edges. There are 7 non-isomorphic tournaments. Altogether there are $1+4+3+7=15$ possibilities.

In case $t(T_7)=10$ the expected number of directed triangles in a random Steiner triple
system is precisely 2. However, the distribution is easily seen to be non-constant
(e.g., the variance is positive). Thus, there is a Steiner triple system with less than two directed triangles.
Namely, $P(T_7) \geq 6$ in this case. \endpf

Fix $T_n$, and let $3 \leq m \leq n$.
Let $T_m$ be a randomly chosen $m$-vertex induced subgraph of $T_n$.
Let $X=a(T_m)$ denote the random variable corresponding to the
number of transitive triples of $T_m$, and let $E[X]$ denote the expectation of $X$.
\begin{prop}
\label{p23}
$E[X] \geq \frac{3}{4}\frac{n-3}{n-2}{m \choose 3}$.
\end{prop}
{\bf Proof}\,
A specific triple of $T_n$ belongs to precisely ${{n-3} \choose {m-3}}$ induced subgraphs
on $m$ vertices. Thus, by (\ref{e1}),
$$
E[X] =\frac{a(T_n) { {n-3} \choose {m-3}}  }{ {n \choose m}}= a(T_n)\frac{m(m-1)(m-2)}{n(n-1)(n-2)}
\geq \frac{3}{4}\frac{n-3}{n-2}{m \choose 3}.
$$
\endpf

{\bf Proof of Theorem \ref{t1}:}\,
Let $n > N+41$ where $N$ is the constant from lemma \ref{l21}.
Let $T_n$ be a fixed $n$-vertex tournament.
We may assume that $n \equiv 1 \bmod 42$, since otherwise we may
delete at most $41$ vertices, apply the theorem on the smaller graph,
and this will not affect the claimed asymptotic number of transitive triples in
the original graph. By Proposition \ref{p23}, the expected number of transitive triples in a random $T_7$
of $T_n$ is at least $26.25(n-3)/(n-2)=26.25(1-o_n(1))$. Hence, the expected number of directed triangles
is at most $8.75(1+o_n(1))$.

Let $p_1$ denote the probability that a random $T_7$ has $t(T_7) \leq 4$. Let $p_2$ denote the
probability that a random $T_7$ has $5 \leq t(T_7) \leq 11$. Let$p_3$ denote the probability that
a random $T_7$ has $t(T_7) \geq 12$. Clearly, $p_1+p_2+p3=1$ and
$$
5p_2+12p_3 \leq 8.75(1+o_n(1)).
$$
Let $Y$ denote the random variable corresponding to $P(T_7)$.
By definition of $p_1,p_2,p_3$ and by Lemma \ref{l22} we have
$$
E[Y]  \geq 7p_1+6p_2+5p_3.
$$
Minimizing $E[Y]$ subject to $p_1+p_2+p_3=1$, $p_i \geq 0$ and $5p_2+12p_3 \leq 8.75(1+o_n(1))$.
yields $p_1=0$, $p_2=13/28(1-o_n(1))$, $p_3=15/28(1+o_n(1))$ and $E[Y] \geq \frac{153}{28}(1-o_n(1))$.

Let $S$ be a fixed $K_7$-decomposition of $K_n$ into ${n \choose 2}/21$ edge-disjoint copies of $K_7$.
By Lemma \ref{l21} such an $S$ exists. Each $s \in S$ corresponds to a $7-set$ of $[n]$.
Let $\sigma$ be a random permutation of $[n]$ and let $S_\sigma$ denote the $T_7$-decomposition
of $T_n$ corresponding to $S$ and $\sigma$. Namely, for each $s \in S$ the corresponding $T_7$-subgraph
of $T_n$, denoted $s_\sigma$, consists of the 7 vertices $\{\sigma(i) ~: ~ i \in s\}$.
Notice that since $\sigma$ is a random permutation, $s_\sigma$ is a random $T_7$ of $T_n$.
Thus, the expected number of edge-disjoint transitive triples of $s_\sigma$ is at least $\frac{153}{28}(1-o_n(1))$.
By linearity of expectation we get that
$$
P(T_n) \geq \frac{{n \choose 2}}{21}\frac{153}{28}(1-o_n(1))=\frac{51}{392}n^2(1+o_n(1)) > 0.13n^2(1+o_n(1)).
$$
\endpf

\section{Edge-disjoint transitive triples in a random tournament}
A random tournament with $n$ vertices is obtained by selecting the orientation of each edge
by flipping an unbiased coin, where all ${n \choose 2}$ choices are independent.
Assume, therefore, that $T_n$ is a random tournament.
\begin{prop}
\label{p31}
$$
{\rm Prob} \left[P_k(T_n) \geq \frac{1}{k(k-1)}n^2(1-o_n(1)) \right] \geq 1-o_n(1).
$$
\end{prop}
{\bf Proof}\,
Let $(x,y)$ be any edge of $T_n$. Clearly, each $K_k$ containing $(x,y)$
induces a $TT_k$ with probability $k!/2^{k \choose 2}$.
Hence, letting $n(x,y)$ denote the number of
transitive $k$-vertex tournaments containing $(x,y)$, we have $E[n(x,y)]={{n-2} \choose {k-2}}k!/2^{k \choose 2}$.
As any two $k$-vertex tournaments containing $(x,y)$ share at most $k-3$ vertices (other than $x$ and $y$) there is
limited dependence between the tournaments containing $(x,y)$ (in fact, for $k=3$ there is complete independence).
Hence, standard large deviation arguments for limited dependence yield that for every $0.5 > \epsilon > 0$,
$$
{\rm Prob}\left[ \left| n(x,y)-{{n-2} \choose {k-2}}\frac{k!}{2^{k \choose 2}}\right| > n^{k-2-\epsilon}\right] =o(n^{-2}).
$$
Thus, with probability $1-o_n(1)$, all edges of $T_n$ lie on at least
$(k(k-1)/2^{k \choose 2})n^{k-2}(1-o_n(1))$
copies of $TT_k$ and at most $(k(k-1)/2^{k \choose 2})n^{k-2}(1+o_n(1))$
copies of $TT_k$.

Consider the ${k \choose 2}$-uniform hypergraph $H$ whose $N={n \choose 2}$ vertices are the edges of $T_n$
and whose edges are the (edge sets of) $TT_k$ copies of $T_n$.
The degree of all the vertices in this hypergraph
is $(k(k-1)/2^{k \choose 2})n^{k-2}(1 \pm o_n(1))=2^{k/2-1-k(k-1)/2}k(k-1)N^{k/2-1}(1 \pm o_n(1))$,
(i.e. the hypergraph is almost regular).
Furthermore, the co-degree of any two vertices in this hypergraph is at most $O(n^{k-3})=O(N^{k/2-1.5})=o(N^{k/2-1})$.
By the result of Frankl and R\"odl \cite{FrRo},
this hypergraph has a matching that covers all but at most $N(1-o_N(1))$ vertices.
Such a matching corresponds to a set of $\frac{1}{k(k-1)}n^2(1-o_n(1))$ edge-disjoint
copies of $TT_k$ in $T_n$. \endpf

\section{Concluding remarks}
\begin{itemize}
\item
Whenever $P(T_n) = n(n-1)/6$ we say that $T_n$ has a {\em transitive} Steiner
triple system. Clearly, this may occur only if $K_n$ has a Steiner triple system,
namely, when $n \equiv 1,3 \bmod 6$. It would be interesting to characterize the
tournaments that have a transitive Steiner triple system.
\item
Conjecture \ref{c1}, if true, would be best possible. We show $f(n) \leq \lceil n(n-1)/6-n/3 \rceil$.
Let $T_3(n)$ be the complete $3$-partite Tur\'an graph with $n$ vertices. It is well-known that
$T_3(n)$ has ${n \choose 2}-\lceil n(n-1)/6-n/3 \rceil$ edges. Denote the partite classes
by $V_1,V_2,V_3$. Orient all edges between $V_1$ and $V_2$ from $V_1$ to $V_2$.
Orient all edges between $V_2$ and $V_3$ from $V_1$ to $V_2$.
Orient all edges between $V_1$ and $V_3$ from $V_3$ to $V_1$.
Complete this oriented graph to a tournament $T_n$ by adding directed edges between any
two vertices in the same partite class in any arbitrary way. Notice that each transitive triple
in $T_n$ contains at least one edge with both endpoints in the same vertex class.
Hence, $P(T_n) \leq \lceil n(n-1)/6-n/3 \rceil$.
\item
Conjecture \ref{c1} has been verified for $n \leq 8$. The values $f(1)=f(2)=f(3)=0$
and $f(4)=1$ are trivial. The values $f(5) = 2$, $f(6)=3$ are easy exercises.
The value $f(7) \geq 5$ is a consequence of Lemma \ref{l22}, and thus $f(7) = 5$ by the above Tura\'n
graph argument. The value $f(8) \geq 7$ (and hence $f(8)=7$) is computer verified.
\item
Conjecture \ref{c1} claims, in particular, that one can cover almost all edges of $T_n$ with
edge-disjoint transitive triples. Proposition \ref{p31} asserts that this is true for the random tournament
and that, in fact, the random tournament can be covered almost completely with edge-disjoint copies of $TT_k$
for every fixed $k$. However, for $k \geq 4$ there are constructions showing that a significant amount
of edges must be uncovered by any set of edge-disjoint $TT_k$. Consider $TT_4$. It is well-known
(cf. \cite{RePa}) that there is a unique $T_7$ with no $TT_4$. Consider the complete $7$-partite digraph
with $n$ vertices obtained by blowing up each vertex of this unique $T_7$ with $n/7$ vertices.
Add arbitrary directed edges connecting two vertices in the same vertex class to obtain a $T_n$.
Clearly, any $TT_4$ of this $T_n$ must contain an edge with both endpoints in the same vertex class.
Hence, $f_4(n) \leq P_4(T_n) \leq 7{{n/7} \choose 2}=O(\frac{1}{14}n^2)$.
Hence at least ${n \choose 2} - 6 f_4(T_n) \geq \frac{1}{14}n^2(1+o(1))$ must be uncovered.
Similar constructions exist for all $k \geq 4$, where the fraction of covered edges tends to
zero as $k$ increases.
\item
It is possible to slightly improve the constant appearing in Theorem \ref{t1}.
Recall that the proof of Theorem \ref{t1} assumed a worst case of $p_1 \geq 0$, where $p_1$ is
the probability that a random $T_7$ has at most four directed triangles.
However, it is very easy to prove that for $n$ sufficiently large, $p_1 > c > 0$ where $c$ is some (small) absolute
constant. This follows from the fact that every $T_{54}$ contains a $TT_7$ \cite{Sa}.
Thus there exists a positive constant $c'$ such that for $n$ sufficiently large,
$T_n$ has at least $c'n^7$ copies of $TT_7$. Hence, a random induced $7$-vertex subgraph of $T_n$ is
a $TT_7$ with constant positive probability. This improvement for $p_1$ immediately implies a
(very small) improvement for the constant appearing in Theorem \ref{t1}.
\end{itemize}

\section*{Acknowledgment}
The author thanks Noga Alon for useful comments.


\begin{thebibliography}{99}

\bibitem{Bo} B. Bollob\'as,
{\em Extremal Graph Theory}, Academic Press, 1978.

\bibitem{FrRo} P. Frankl and V. R\"odl,
{\em Near perfect coverings in graphs and hypergraphs},
European J. Combinatorics 6 (1985), 317-326.

\bibitem{Ga} R.B. Gardner,
{\em Optimal packings and coverings of the complete directed graph with $3$-circuits and with transitive triples},
Congressus Numerantium 127 (1997), 161-170.

\bibitem{PhLi}K.T. Phelps and C.C. Lindner,
{\em On the number of Mendelsohn and transitive triple systems},
European. J. Combin. 5 (1984), 239-242.

\bibitem{RePa} K.B. Reid and E.T. Parker,
{\em Disproof of a conjecture of Erd\H{o}s and Moser on tournaments},
J. Comb. Theory 9 (1970), 225-238.

\bibitem{Sa} A. Sanchez-Flores,
{\em  On tournaments free of large transitive subtournaments},
Graphs and Combinatorics 14 (1998), 181--200.

\bibitem{Sk} D.B. Skillicorn, 
{\em A note on directed packings of pairs into triples},
Ars Combin. 13 (1982), 227-229. 

\bibitem{Wi} R. M. Wilson,
{\em Decomposition of complete graphs into subgraphs isomorphic to a given
graph}, Congressus Numerantium XV (1975), 647-659.

\end{thebibliography}
\end{document}